# Traffic flow brake light model simulation based on driver behavior learning

**Rui Shen[1, *], Siqi Yuan[2], Zhihua Kong[2]**

[1]School of Beijing Jiaotong University, Beijing, China
[2]School of aaaa, bbbb University, Shanghai, China

*Corresponding author: zzzz@xxxx.com

**Abstract**. The theory of urban traffic flow has been developed and new types of meta-automata have emerged and simulate realistic traffic conditions relatively well. Among these models, the brake light model can simulate the three-phase traffic flow theory very well. However, the existing brake light model also has certain shortcomings, in that the model will change the speed of congestion propagation upward when the model is covariant, which is not realistic, and the model also lacks simulation parameters for driver behavior. In this paper, we propose a new model based on the brake light model, which can achieve a certain degree of simulation of driver behavior by adjusting the parameters, and also achieve a stabilization of the propagation speed of the congestion wave when the parameters are changed.

## 1. Introduction

In recent years, China's urban economy has been developing and the urbanization process has been accelerating, in which the development of traffic has become increasingly important. Scholars at home and abroad have carried out simulation for various characteristics of traffic flow, and proposed the use of metacellular automata to simulate traffic flow from the microscopic level, which has achieved obvious results. However, the microscopic characteristics of traffic flow still have many phenomena that are difficult to be simulated by ordinary meta-cellular automata models, such as the back-lag phenomenon and traffic collapse phenomenon. At the same time, there are also phenomena that are difficult to match the actual driver's personality when traffic flow is simulated. The KKW and KKS models of Kerner[1] et al. add the assumption that the vehicle will tend to maintain synchronous speed with the vehicle in front while ensuring safety, while Davis[2] proposes to add variables to the vehicle-following model to simulate the traffic collapse phenomenon. Chmura[3] et al. propose to limit the acceleration and deceleration of vehicles in the model to a reasonable range and successfully simulate a variety of important synchronous flow characteristics of traffic flow. Tian et al[4] , on the other hand, used the brake light element in the model and proposed a brake light model, which was able to simulate phenomena such as free flow, synchronous flow, and wide motion blockage.

In general, most of the current research has focused on the use of adding reasonable parameters to make the traffic flow simulation simulation more close to reality, but there is a lack of simulation of driver behavior, especially there is less research related to vehicle driving behavior. Therefore, this paper adopts further optimization of the DTGBLM model by adding parameters that reflect the driver in order to study the adoption of driving behavior to achieve a more reasonable simulation of the actual traffic flow.

**2. Another section of your paper**
Nagel and Schreckberg in 1992 proposed a meta-cellular automaton model for simulating vehicular traffic, the NaSch model (as it has also been called), by means of a simple rule to simulate a follow-the-road model of vehicles. Time, space and vehicle speed are discretized by integers. Roads are divided into discrete lattices, or tuplets, of equal distance. Each tuple is either empty or occupied by a vehicle. The speed of the vehicle can take values between (0 and Vmax). The model evolves according to the following rules

The update is performed during the process from moment t to moment t+1 according to the following rule.

(1) Random slowing.

$$V_n(t+1) = \begin{cases} \max(V_{n(t+1)} - b, 0), & \text{if } r < p \\ V_{n(t+1)}, & \text{otherwise} \end{cases}$$

Slowing down with random probability p.

Rule (1) introduces stochastic slowing to reflect differences in driver behavior, so that both stochastic acceleration behavior and overreaction behavior during deceleration can be captured. This rule is also a crucial factor in the spontaneous generation of blockages.

(2) Acceleration.
$$V_n(t+1) = \min(V_n(t) + a, V_{max})$$
Rule (2) reflects the driver's tendency to drive at the maximum possible speed.

(3) Deceleration.
$$V_n(t+1) = \min(d_n(t) - c, 0)$$
Rule (3) Ensure that the vehicle does not collide with the vehicle in front of it.

(4) Location updates.
$$X_{n(t+1)} = X_{n(t)} + V_{n(t+1)}$$
The vehicle moves forward at the updated speed.

Of which $X_n$, $V_n$ denote the nth vehicle position and speed, respectively; b (b ⩾ 1) is the vehicle length; d denotes the number of empty metacells between vehicle n and the preceding vehicle n + 1; p denotes the random slowing probability; and $V_{max}$ is the maximum velocity.

When receiving the paper, we assume that the corresponding authors grant us the copyright to use the paper for the book or journal in question. Should authors use tables or figures from other Publications, they must ask the corresponding publishers to grant them the right to publish this material in their paper.

**3. Introduction of the DTGBLM model and the new model**
The brake light model adds a completely new zero-one variable to the base NaSch model $S_n(t)$ represents the state of the nth vehicle's brake light on at time t. 1 represents the vehicle's brake, etc., coming on and the rear vehicle noticing the deceleration, and 0 represents the vehicle not decelerating. Its model update rule is as follows.

Step 1: Determine the random slowing probability

$$p_n = p(v_n(t), s_{n+1}(t), t_{n,h}(t), t_{n,sa}(t)) = \begin{cases} p_b: \text{If } s_{n+1} = 1 \text{ and } t_{n,h}(t) < t_{n,sa}(t) \\ p_0: \text{If } v_n(t) = 0 \\ p_d: \text{otherwise} \end{cases} \quad (1)$$

$$s_{n+1}(t+1) = 0 \qquad (2)$$

Step 2: Acceleration Step

If $s_{n+1}(t) = 0$ or $t_{n,h}(t) \geq t_{n,sa}(t)$, and satisfies $v_n(t) = 0$, then $v_n(t+1) = \min(v_n(t) + a_1, v_{max}, \lceil d_{n,eff}(t)/T \rceil)$; otherwise $v_n(t+1) = \min(v_n(t) + a_2, v_{max}, \lceil d_{n,eff}(t)/T \rceil)$

If $v_n(t+1) < v_n(t)$, then the brake lights of vehicle n come on at moment $t+1$, i.e. $s_n(t+1) = 1$

Step 3: Random Slowing

Given a random number rand() within [0,1], if $\text{rand}() < p_n$, then $v_n(t+1) = \max(v_n(t+1) - b_{rand}, 0)$; if $p_n = p_b$, then $s_n(t+1) = 1$

Step 4: Location update

$$x_n(t+1) = x_n(t) + v_n(t+1) \qquad (3)$$

where $p_n$ is the random slowing probability of vehicle n at moment t+1, and $d_n(t) = x_{n+1} - x_n - L_{veh}$ is the spatial distance between vehicle n and vehicle n+1, and $s_n(t)$ denotes the state of the brake lights of vehicle n at moment t ($s_n(t) = 0$ indicates that the brake lights are off and $s_n(t) = 1$ denotes the brake light is on), and $t_{n,h}(t) = d_n(t)/v_n(t)$ denotes the time distance between vehicle n and the vehicle n+1 in front of it at moment t, and $t_{n,sa}(t)$ denotes the safe time distance between vehicle n and the preceding vehicle n+1 at moment t, and h denotes the maximum time distance for the preceding vehicle to have an effect on the following vehicle, thus $t_{n,sa}(t) = \min(v_n(t), h)$. $d_{n,eff}(t) = d_n(t) + \max(v_{anti} - b_{anti}, 0)$ is the effective spatial distance between vehicle n and vehicle n+1, and $v_{anti} = \min(d_{n+1}, v_{n+1})$ is the desired velocity of vehicle n+1, and $b_{anti}$ is the desired deceleration of the front vehicle. Here the acceleration and deceleration process of the front vehicle is considered on the basis of the NaSch model, and the spatial distance between vehicle n and vehicle n+1 $d_n(t)$ is redefined as $d_{n,eff}(t)$. T is the expected headway time distance ($T > 1s$), and $p_b, p_0, p_d$ is the random slowing probability for the different cases of $a_1$ and $a_2$ are two acceleration values used to ensure that the simulated synchronous traffic flow state matches the actual one. Only when $b_{anti} \geq b_{rand}$ time the accident can be avoided.

Step 1: Determine driver characteristics

$$I_n = \begin{cases} 1 \\ 0 \end{cases} \qquad (4)$$

$I_n$ represents the driving habits of the driver, and when $I_n = 0$, the driver is driving normally, and when $I_n = 1$, the driver highlights the characteristic of impatience or fear, causing a greater change in acceleration or deceleration.

Step 2: Determine the random slowing probability P

$$p_n = p(v_n(t), s_{n+1}(t), t_{n,h}(t), t_{n,sa}(t)) = \begin{cases} p_b: \text{If } s_{n+1} = 1 \text{ and } t_{n,h}(t) < t_{n,sa}(t) \\ p_0: \text{If } v_n(t) = 0 \\ p_{d1}: \text{Otherwise and } I_n = 1 \\ p_{d2}: \text{Otherwise and } I_n = 0 \end{cases} \qquad (5)$$

$$s_{n+1}(t+1) = 0 \qquad (6)$$

Step 3: Acceleration Step

If $s_{n+1}(t) = 0$ or $t_{n,h}(t) \geq t_{n,sa}(t)$, the

then

$$\begin{cases} v_n(t+1) = \min(v_n(t) + a_1, v_{max}, \lceil d_{n,eff}(t)/T \rceil): & if\ I_n = 1 \\ v_n(t+1) = \min(v_n(t) + a_2, v_{max}, \lceil d_{n,eff}(t)/T \rceil): & if\ I_n = 0 \end{cases} \quad (7)$$

Otherwise..

$$v_n(t+1) = \min(v_n(t) + a_3, v_{max}, \lceil d_{n,eff}(t)/T \rceil) \quad (8)$$

If $v_n(t+1) < v_n(t)$, then the brake lights of vehicle n come on at moment t + 1, i.e. $s_n(t+1) = 1$

Step 4: Random Slowing
Given a random number rand() within [0,1], if $rand() < p_n$, the
(a) The rule. $v_n(t+1) = \max(v_n(t+1) - b_{rand}, 0)$.
If $p_n = p_b$, then $s_n(t+1) = 1$

Step 5: Location update

$$x_n(t+1) = x_n(t) + v_n(t+1) \quad (9)$$

## 4. Model simulation

### 4.1. DTGBLM model simulation

Referring to some of the literature, finally this paper identifies a metacell model with periodic boundaries and the road length is taken to be 2500 cells and the model specific model parameters are shown in the following table.

**Table 1.** Three Scheme comparing.

| Parameters | $L_{cell}$ | $L_{veh}$ | $V_{max}$ | h | T | $p_b$ | $p_o$ | $p_d$ | a | $b_{rand}$ | $b_m$ | $v_{cri}$ |
|---|---|---|---|---|---|---|---|---|---|---|---|---|
| Units | m | $L_{cell}$ | $L_{cell}/s$ | s | s | - | - | - | $L_{cell}/s^2$ | $L_{cell}/s^2$ | $L_{cell}/s^2$ | $L_{cell}/s$ |
| Value | 1.5 | 5 | 20 | 6 | 2.8 | 0.8 | 0.45 | 0.01 | 1 | 1 | 5 | 5 |

According to the above parameters, the random slowing rate of vehicles by adjusting $p_o$ and road vehicle density to achieve the detection of road traffic. The DTGBLM model as well as the new model is programmed and numerically simulated according to the previous presentation.

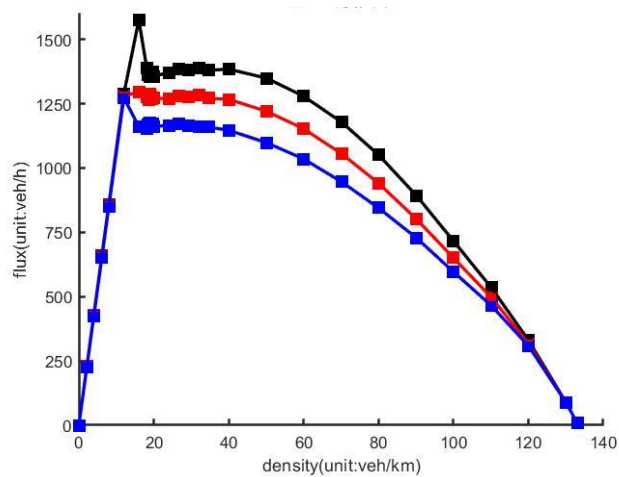

**Figure 1.** Two or more references.

From the figure it can be seen that at different slowing rates $p_d$ the DTGBLM model shows not much difference, after the traffic flow density after 27veh/km, the slowing rate $p_d$ places a limit on the upper limit of the curve, after which the traffic flow shows a collapse after 36veh/km, which is a good representation of the three-phase traffic flow theory of traffic flow.

The following figure shows the trajectory roadmap obtained by varying the slowing rate $p_d$ The trajectory roadmap obtained, the figure takes the last 500 time steps of the simulation 10000 time steps, it can be seen that different slowing rates have a great impact on the vehicle propagation speed, which is not in line with the actual situation.

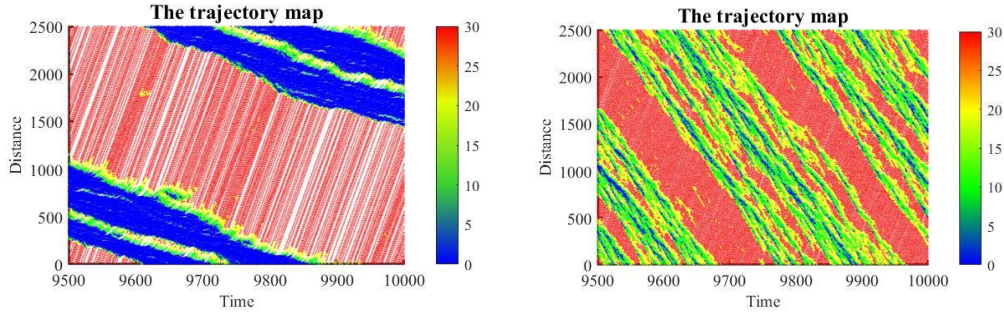

**Figure 2.** Two or more references.

*4.2. The new model simulation*
By referring to some of the literature, it was confirmed that the new model road length was taken to be the same 2500 element cells and the model specific model parameters are shown in the following table.

**Table 2.** Three Scheme comparing.

| Parameters | $L_{cell}$ | $L_{veh}$ | $V_{max}$ | h | T | $p_b$ | $p_o$ | $p_d$ | a | $b_{rand}$ | $b_m$ | $v_{cri}$ |
|---|---|---|---|---|---|---|---|---|---|---|---|---|
| Units | m | $L_{cell}$ | $L_{cell}/s$ | s | s | - | - | - | $L_{cell}/s^2$ | $L_{cell}/s^2$ | $L_{cell}/s^2$ | $L_{cell}/s$ |
| Value | 1.5 | 5 | 20 | 6 | 2.8 | 0.8 | 0.45 | 0.01 | 1 | 1 | 5 | 5 |

The new model flow versus density is shown in Fig. From the figure, we can see that the new model is able to simulate the change from free flow to synchronous flow by the image of the three-phase traffic flow achieved at different slowing rates for the same flow density. The following basic and spatio-temporal diagrams are obtained.

*4.2.1. Flow density analysis.* The following figure depicts the relationship between traffic flow and density, from which it can be seen that the traffic flow gradually collapses in the second phase as the slowing rate increases, proving that the new model is able to simulate not only the free flow to synchronous flow but also the sudden change of flow from synchronous flow to congested flow. Moreover, under the same conditions, the model gradually shows flow differences only after the density exceeds 36veh/km, proving that the individualized distribution of drivers facilitates, within certain limits, the smooth operation of the traffic flow.

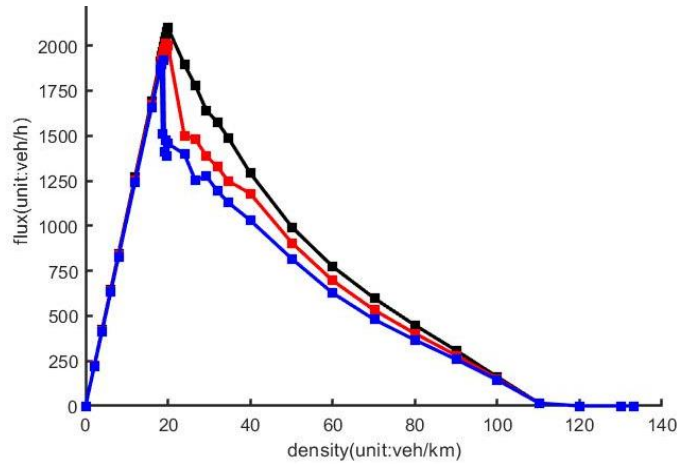

**Figure 3.** Two or more references.

*4.2.2. The trajectory map analysis.* Similarly, the trajectory roadmap of the new model depicts, for time steps starting at 9500 and ending at 10000, a total of 500 time steps in space and time, at different slowing rates, in the figure from which the spontaneously generated blockage phenomenon can be simulated and a wide motion blockage occurs, which is a more rea onable simulation of the actual road.

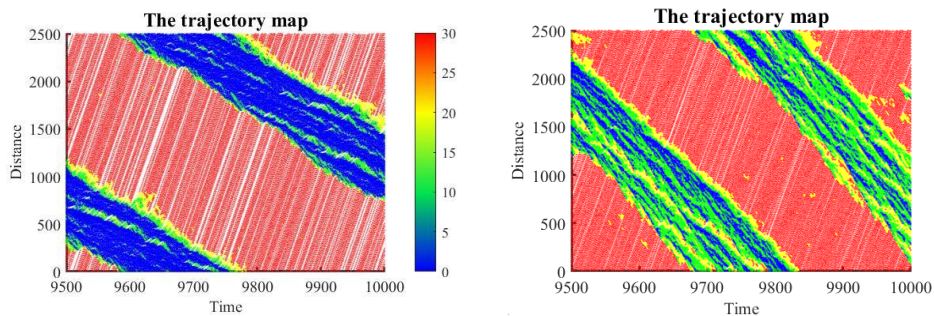

**Figure 4.** Two or more references.

## 5. Conclusions

In this paper, by adding the study of driver characteristics to the original brake light metacell model, not only, retaining the characteristics of the brake light metacell automaton capable of simulating free flow, synchronous flow and wide motion blockage, but also, in the case of comparison, further understanding the impact of each parameter in the specific model on the model, achieving an in-depth and comprehensive understanding of the new role brought by the new variables in the establishment of the new model. The following conclusions are drawn.

(1) The brake light model can simulate the basic state of some traffic flows better, but in the face of complex traffic systems, the brake light model is difficult to achieve the simulation of congestion propagation of actual traffic flows.

(2) The individualized distribution of vehicle drivers at different slowing rates facilitates the smooth operation of traffic flow and promotes better traffic flow within certain limits.

(3) Traffic collapse phenomenon is more difficult to occur under different driving conditions. When there are fewer random vehicle decelerations, the traffic flow will gradually decline after reaching its peak and no traffic collapse will occur, while when there are more random vehicle decelerations, the traffic flow is highly susceptible to traffic collapse and traffic efficiency is easily affected.

**References**
[1] KERNER B S. Three-Phase Theory of City Traffic: Moving Synchronized Flow Patterns in under-Saturated City Traffic at Signals[J]. Physica A: Statistical Mechanics and Its Applications, 2014, 397: 76–110. DOI:10.1016/j.physa.2013.11.009.
[2] DAVIS L C. Multilane simulations of traffic phases[J]. Physical Review E - Statistical, Nonlinear, and Soft Matter Physics, 2004, 69(1 2): 161081–161086.
[3] CHMURA T, HERZ B, KNORR F, 等. A simple stochastic cellular automaton for synchronized traffic flow[J]. Physica A: Statistical Mechanics and its Applications, 2014, 405: 332–337. DOI:10.1016/j.physa.2014.03.044.
[4] TIAN J, TREIBER M, MA S, 等. Microscopic driving theory with oscillatory congested states: Model and empirical verification[J]. Transportation Research Part B: Methodological, 2015, 71: 138–157. DOI:10.1016/j.trb.2014.11.003.
[5] LI Q-L, JIANG R, DING Z-J, 等. A new cellular automata traffic flow model considering asynchronous update of vehicle velocity[J]. International Journal of Modern Physics C, 2020, 31(12). DOI:10.1142/S0129183120501673.